\documentclass[12pt]{article}

\usepackage{epsfig}
\usepackage{amssymb}
\usepackage{amsmath}
\usepackage{subfigure}

\def\cond{\mbox{cond\/}}

\def\Cst{{\mathbb C}}

\def\Nst{{\mathbb N}}

\def\Rst{{\mathbb R}}
\def\Rdst{{\Rst^d}}

\def\Tst{{\mathbb T}}
\def\Zst{{\mathbb Z}}

\def\Csp{{\boldsymbol C}}

\def\Hsp{{\boldsymbol H}}

\def\lsp{{\boldsymbol\ell}}
\def\lisp{{\lsp^1}}
\def\liZ{{\lisp(\Zst)}}

\def\ltsp{{\lsp^2}}
\def\ltZ{{\lsp^2(\Zst})}

\def\Lsp{{\boldsymbol L}}

\def\Ltsp{{\Lsp^2}}
\def\LtR{{\Ltsp(\Rst)}}

\def\Lbsp{{\boldsymbol{\mathcal L}}}

\def\Xsp{{\boldsymbol X}}

\def\Ysp{{\boldsymbol Y}}

\def\slb{{{\raise 0.5pt \hbox{\footnotesize $[$}}}}
\def\srb{{{\raise 0.5pt \hbox{\footnotesize $]$}}}}

\def\slp{{{\raise 0.5pt \hbox{\footnotesize $($}}}}
\def\srp{{{\raise 0.5pt \hbox{\footnotesize $)$}}}}

\def\slcb{{{\raise 0.5pt \hbox{\footnotesize $\{$}}}}
\def\srcb{{{\raise 0.5pt \hbox{\footnotesize $\}$}}}}

\def\sabs{{{\raise 0.5pt \hbox{\footnotesize $|$}}}}

\def\CHI{\hbox{\raise .5ex \hbox{$\chi$}}}

\def\shadowbox{\hbox{\rule[-0.0ex]{0.1ex}{1.2ex}%
\hspace{-0.1ex}\rule[-0.0ex]{1.2ex}{0.1ex}%
\hspace{0.0ex}\rule[-0.0ex]{0.1ex}{1.2ex}\hspace{-1.3ex}%
\rule[1.15ex]{1.25ex}{0.1ex}\hspace{-0.0ex}\rule[-0.25ex]{0.3ex}{1.1ex}%
\hspace{-1.2ex}\rule[-0.25ex]{1.1ex}{0.25ex}}}

\def\shadowbox{\hbox{\rule[-0.0ex]{0.1ex}{1.2ex}%
\hspace{-0.1ex}\rule[-0.0ex]{1.2ex}{0.1ex}%
\hspace{0.0ex}\rule[-0.0ex]{0.1ex}{1.2ex}\hspace{-1.3ex}%
\rule[1.15ex]{1.25ex}{0.1ex}\hspace{-0.0ex}\rule[-0.25ex]{0.3ex}{1.1ex}%
\hspace{-1.2ex}\rule[-0.25ex]{1.1ex}{0.25ex}}}
\def\qed{\ifmmode \hbox{\hfill\shadowbox}
     \else \hphantom{x}\hfill\shadowbox \fi}

\newtheorem{theorem}{Theorem}[section]

\newtheorem{proposition}[theorem]{Proposition}
\newtheorem{corollary}[theorem]{Corollary}

\newcommand{\lam}{\lambda}
\newcommand{\lamkl}{\lambda^{|k-l|}}
\newcommand{\gmn}{g_{m,n}}
\newcommand{\gnm}{g_{m}^{(N)}}
\newcommand{\gm}{g_{m}}
\newcommand{\gamn}{\gamma_{m,n}}

\newcommand{\cgm}{{^{\footnotesize P}\!\!g_{m}}}
\newcommand{\cgmn}{{^{\footnotesize P}\!\!g_{m,n}}}
\newcommand{\gam}{\gamma_{m}}
\newcommand{\go}{\gamma}

\newcommand{\gnam}{\gamma_{m}^{(N)}}
\newcommand{\gnamn}{\gamma_{m,n}^{(N)}}
\newcommand{\cnam}{{^{\footnotesize P}\!\!\gamma_{m}^{(N)}}}

\newcommand{\en}{e^{(N)}_m}
\newcommand{\phin}{\phi^{(N)}}
\newcommand{\phinm}{\phi^{(N)}_m}
\newcommand{\phim}{\phi_m}
\newcommand{\phimn}{\phi_{m,n}}
\newcommand{\EN}{E_N}
\newcommand{\SN}{S_N}
\newcommand{\SI}{S^{-1}}
\newcommand{\SR}{S^{\frac{1}{2}}}
\newcommand{\SRI}{S^{-\frac{1}{2}}}

\newcommand{\SNI}{S^{-1}_N}
\newcommand{\CN}{{^{\footnotesize P}\!S_N}}
\newcommand{\CNI}{{^{\footnotesize P}\!S_{N}^{-1}}}
\newcommand{\PN}{P_N}
\newcommand{\PNg}{P_N g}
\newcommand{\skap}{\sqrt{\kappa}}
\newcommand{\HN}{\Hsp_{N}}
\newcommand{\HNN}{\Hsp_{N \times N}}
\newcommand{\kap}{\kappa}
\newcommand{\lra}{\Leftrightarrow}
\newcommand{\xx}{\times}
\newcommand{\aaa}{\tau}
\newcommand{\AAA}{T}
\newcommand{\DDD}{a}

\begin{document}

\title{Rates of convergence for the approximation of dual 
shift-invariant systems in $\ltZ$}
%and wavelet frames}
%\title{\bf Rates of convergence for the approximation of dual Gabor 
%and wavelet frames}
%\title{\bf Rates of convergence for the approximation of dual Gabor and
%wavelet frames using finite-dimensional models}
%\title{\bf Gabor frames, filter banks, and the projection method}
%\title{Gabor frames, the finite section method and approximation}
\author{Thomas Strohmer\thanks{Department of Mathematics, University
of California, Davis, CA-95616, email: strohmer@math.ucdavis.edu.}}

\date{}
\maketitle

\begin{abstract}
A shift-invariant system is a collection of functions $\{\gmn\}$ of the form
$\gmn(k) = \gm(k-an)$. Such systems play an important
role in time-frequency analysis and digital signal processing. A principal
problem is to find a dual system $\gamn(k) = \gam(k-an)$ such that each 
function $f$ can be written as $f = \sum \langle f, \gamn \rangle \gmn$.
The mathematical theory usually addresses this problem in infinite dimensions
(typically in $\LtR$ or $\ltZ$), whereas numerical methods have to operate
with a finite-dimensional model. Exploiting the link between the frame 
operator and Laurent operators with matrix-valued symbol, we apply the
finite section method to show that the 
dual functions obtained by solving a finite-dimensional problem converge 
to the dual functions of the original infinite-dimensional problem in $\ltZ$. 
For compactly supported $\gmn$ (FIR filter banks) we prove an exponential 
rate of convergence and derive explicit expressions for the involved 
constants. Further we investigate under which conditions one can replace 
the discrete model of the finite section method by the periodic 
discrete model, which is used in many numerical procedures. Again we 
provide explicit estimates for the speed of convergence. Some remarks 
on tight frames complete the paper.
\end{abstract}
\mbox{Key words:} Shift-invariant systems, finite section method, 
block Toeplitz matrices, Laurent operator, Gabor frame, filter banks,
band matrices.

\section{Introduction}
\label{s:intro}

Shift-invariant systems play an important role in time-frequency
analysis and digital signal processing. In this
paper we consider time-discrete shift-invariant systems, by which we mean 
functions $\gmn$ of the form
\begin{equation}
\gmn(k) = \gm(k - n \DDD), \qquad k,n \in \Zst, m = 0,1,\dots,M-1
\label{sis}
\end{equation}
where $\gm \in \ltZ$ and $0 < \DDD \in \Nst$.
In the engineering literature such a system is known as {\em filter bank}. 

Since the $\gmn$ are in general not orthogonal, we are interested in
finding dual systems 
\begin{equation}
\gamn(k) = \gam(k-n \DDD), \qquad k,n \in \Zst, m = 0,1,\dots,M-1
\label{dualsis}
\end{equation}
such that any $f \in \ltZ$ has a $\ltZ$-convergent representation
\begin{equation}
f = \sum_{m,n} \langle f, \gamn \rangle \gmn \,.
\label{fseries}
\end{equation}

There is a variety of research papers dealing with shift-invariant
systems in $\LtR$ or $\ltZ$, there is also no lack of methods claiming to
efficiently solve for the duals, those method naturally operate in
finite dimensions. But only little is known about the relation of the 
finite-dimensional solutions and the original solutions of the infinite 
dimensional model. For practical purposes it is not enough to ensure that an
approximation converges to the original solution, but one also should try
to provide a priori bounds on the rate of convergence. 
The goal of this paper is to clarify the relation between the approximate 
dual functions computed by a finite-dimensional approach and the original 
dual functions in $\ltZ$. 

Throughout the paper we will use two key ideas to derive our results.
We exploit the link between frame operators of shift-invariant systems
and Laurent operators with matrix valued symbol to apply the finite section 
method for the proof of convergence. Furthermore we will discuss the rate 
of convergence making use of results of the off-diagonal entries of the
inverse of certain types of matrices.

A detailed investigation of shift-invariant systems in $\LtR$ can be found
in the papers of Ron and Shen~\cite{RS95a,RS97b}. In~\cite{Jan98}
Janssen analyses duality conditions for shift-invariant systems in
$\LtR$ and $\ltZ$ with emphasis on Gabor systems. For an overview of the 
theory of filter banks the reader may consult~\cite{VK95}.

The concept of frames~\cite{DS52,You80} provides a natural environment to 
study shift-invariant systems~\cite{Dau92,RS95a,Jan95}. We say that the
set $\{\gmn\}$ is a frame for $\ltZ$ if there exist constants $A,B>0$ 
such that
\begin{equation}
\label{framedef}
A \|f\|^2 \le \sum_{m,n} |\langle f, \gmn \rangle|^2 \le B \|f\|^2 \,,
\qquad \forall f \in \ltZ\,.
\end{equation}
If the shift-invariant system $\{\gmn\}$ satisfies condition~\eqref{framedef}
we will also call it {\em shift-invariant frame}. The dual frame is 
given by 
\begin{equation}
\gamn(k) = (S^{-1} \gm)(k - n\DDD)
\label{dualframe}
\end{equation}
where $S$ is the frame operator defined by
\begin{equation}
Sf = \sum_{m,n} \langle f, \gmn \rangle \gmn \,.
\label{frameop}
\end{equation}

Among all dual systems satisfying~\eqref{fseries} the duals defined 
by~\eqref{dualframe} have minimal $\ltsp$-norm. In this paper we restrict
our attention to duals with minimal $\ltsp$-norm, therefore we
will be somewhat sloppy in our terminology and refer to the $\gamn$ 
defined in~\eqref{dualframe} as {\em the} duals.

We define the analysis operator $F$ by
\begin{equation}
F: f \in \ltZ \rightarrow Ff = \{ \langle f, \gmn \rangle\}_{m,n}
\label{frameanal}
\end{equation}
and the synthesis operator, which is just the adjoint operator of $F$, by
\begin{equation}
F^{\ast}: c \in \ltsp(\Zst \times {\cal I}) \rightarrow 
F^{\ast} c = \sum_{m,n} c_{m,n} \gmn 
\label{framesyn}
\end{equation}
where ${\cal I}:=\{0,1,\dots,M-1\}$. The frame operator can be written
as $S = F^{\ast} F$.

If the $\gm$ are constructed from a single function $g$ by modulations, i.e., 
\begin{equation}
\gm(k) = e^{2\pi i mk/M} g(k)
\label{gabor} 
\end{equation}
we obtain the well-known {\em Gabor system} or {\em Weyl-Heisenberg
system}~\cite{Gab46,FS98}. A nice property of Gabor frames is that the dual 
frame is also generated by a single function, we have~\cite{Dau90} 
\begin{equation}
\gamn(k) =  e^{2\pi i mk/M} \gamma(k-na)\qquad \text{with}\,\,\,
\gamma = S^{-1}g\,\,.\notag
%\label{dualgabor}
\end{equation}
In the filter bank community systems of the form~\eqref{gabor} are known 
as (oversampled) DFT filter banks~\cite{BH98}.

For some special choices of $\{\gmn\}$ closed-form solutions for the duals 
are known~\cite{Dau90,FZ95,Jan96} (which does not necessarily mean that 
those closed-form solutions can be cheaply computed in practice).
In general however the numerical computation of the dual frame elements $\gam$ 
involves the solution of the system of equations~\eqref{dualframe} of
biinfinite order, independently of the choice of the basis for the matrix 
representation of $S$ (be it the standard basis in $\ltZ$ or
the so-called polyphase representation used in the theory of filter banks). 
We certainly cannot solve a system of biinfinite order, but we can
solve it approximately, based on a finite-dimensional model.  

Whatever finite-dimensional model we choose to compute approximate duals,
we have to ask if these approximations really converge to the original solution
with increasing dimension and, whenever possible, we have to give
a priori bounds on the rate of convergence.
These are the problems we are going to attack in this paper.

Before we proceed we introduce some notation. Let $x \in \ltZ$, then
$\|x\|$ is the standard $\ltsp$-norm of $x$. We define the spectrum of
an operator $T$ as usual by $\sigma(T)=\{\mu \in \Cst: \det(T-\mu I) = 0\}$.
If $T$ is positive definite and invertible then 
$\|T\| = \max \{\mu: \mu \in \sigma(T)\}$ and in this case
$\cond(T):=\|T\|\|T^{-1}\|$ denotes the condition number of $T$.

\section{Dual frames and the finite section method} 
\label{ss:finite}

Due to the shift-invariance property~\eqref{dualsis} of the $\gamn$ we 
only have to compute approximations to the duals $\gam = \gamma_{m,0}$ for
$m=0,\dots,M-1$, since all other $\gamn$ are just translations of $\gam$. 
Equations~\eqref{dualframe} 
and~\eqref{frameop} tell us that all shifted copies $\gmn$ 
are required for the computation of the duals $\gam$. But if the $\gm$ have
nice decay properties, such that most of their energy is contained in 
some interval around the origin, say, then one may hope that the $\gam$ 
will have the same properties. One could further argue that using only those 
functions $\gmn$ ``living'' in some interval around the origin, $[-N,N]$ say,
should be sufficient to compute a good approximation to $\gam$. 

Let us put these vague ideas in a more precise mathematical setting.
For $N \in \Nst$ define the orthogonal projections $\PN$ by
\begin{equation}
\label{defP}
\PN x = (\dots, 0,0, x_{-N}, x_{-N+1},\dots , x_{N-1}, x_N, 0,0, \dots)
\end{equation}
and identify the image of $\PN$ with the $2N+1$-dimensional space 
$\Cst^{N+1}$. We construct a finite-dimensional approximation $\SN$ of the 
frame operator $S$ via the truncated frame elements $\PN \gmn$ by
\begin{equation}
\label{pnspn}
\SN = (\PN F^{\ast})(F \PN)
\end{equation}
where $F$ is defined in~\eqref{frameanal}.
The $N$-th approximation $\gnam$ to $\gam$ is then given by the  
solution of the finite-dimensional system of equations $\SN \gnam = \PN \gm$.
In other words we truncate {\em all} frame elements $\gmn$ to
the interval $[-N,N]$.

%In other words we use only the information of the $\gmn$ contained in the
%interval $[-N,N]$.%, as illustrated in Figure~\ref{fig:gabsys}.
%\begin{figure}
%\begin{center}
%\subfigure[The functions $\gmn$]{
%\epsfig{file=gabsystem.eps,width=66mm,height=55mm}} 
%\subfigure[The functions $\PN \gmn$]{
%\epsfig{file=gabsysproj.eps,width=66mm,height=55mm}} 
%\caption{For the $N$-th approximation we use only the information contained
%in the interval $[-N,N]$, which is illustrated in the lower image.}
%\label{fig:gabsys}
%\end{center}
%\end{figure}

\noindent
The questions that arise with such an approximation scheme are:\\
(i) Does $\gnam$ converge to $\gam$ for $N \rightarrow \infty$? \\
(ii) If the $\gm$ satisfy certain decay conditions, can we give an estimate
on the rate of convergence? In other words: given $\delta>0$, can we provide 
an a priori estimate for the size of the interval $[-N,N]$, such that
$\|\gam - \gnam\| \le \delta$, where $N$ depends on $\delta$ and the
decay of $\gm$?

In numerical procedures people often use a periodic finite-dimensional
model in order to approximate $\gam$. Instead of the truncated
functions $\PN \gm$ (which could be thought of as padded with
zeros beyond the interval $[-N,N]$), the $\PN \gm$ are extended
periodically across the interval $[-N,N]$. In other words all
computations are done in the ring $\Zst_{2N+1} = \Zst\!\! \mod \! (2N+1)$.
In Section~\ref{s:periodic} we will investigate the behavior of the
periodic approximate dual functions, when the length of the
period approaches infinity.

\subsection{The finite section method and Laurent operators} 
\label{ss:laurent}

The {\em finite section method} or {\em projection method} is a classical 
technique to approximate the solution of infinite Toeplitz-type systems. 
For an introduction and many variations of the theme the reader is 
referred to the books of Gohberg and Fel'dman~\cite{GF74}, B\"ottcher and 
Silbermann~\cite{BS90} and Hagen, Roch, and Silbermann~\cite{HRS95}. 
Since we use the finite section method as a main tool 
throughout the paper, we briefly describe it in the specific context 
we will apply it.

Let $\Tst$ be the complex unit circle and let $\aaa(z)$ be a continuous
function on $\Tst$, $\aaa \in \Csp(\Tst)$. Let $\aaa_{k} \in \Cst$ stand for the 
$k$-th Fourier coefficient of $\aaa$,
$$ 
\aaa_{k} = \int \limits_{0}^{1} \aaa(e^{2\pi i \omega}) 
e^{-2\pi ik \omega}\,d\omega
$$
then the operator $\AAA(\aaa)$ which acts on functions $x \in \ltZ$ by
\begin{equation}
(\AAA x)_k = \sum_{j \in \Zst} \aaa_{k-j} x_j \,,\qquad  k=0,1,\dots \,,
\label{laurentop}
\end{equation}
is called {\em Laurent operator} with {\em generating function} or
{\em symbol} $\aaa$. A Laurent operator may also be characterized as a 
bounded linear operator acting on $\ltZ$ which commutes with the forward
shift operator on $\ltZ$, cf.~Chapter~13.2 in~\cite{GGK93}. 
Since a Laurent operator can be represented 
by a biinfinite Toeplitz matrix, we will often identify $\AAA(\aaa)$ 
with its matrix representation.

Consider the operator equation $\AAA x =y, x,y \in \ltZ$ or, in matrix 
representation 
\begin{gather}
\label{toepsys}
\begin{bmatrix}
\ddots & \vdots   & \vdots   & \vdots  &         \\
\dots  & \aaa_{0}  & \aaa_{1}  & \aaa_{2} & \dots   \\
\dots  & \aaa_{-1} & \aaa_{0}  & \aaa_{1} & \dots   \\
\dots  & \aaa_{-2} & \aaa_{-1} & \aaa_{0} & \dots   \\
       & \vdots   & \vdots   & \vdots  & \ddots  \\
\end{bmatrix}\,\,
\begin{bmatrix}
\vdots \\ x_{-1} \\ x_0 \\ x_1 \\ \vdots
\end{bmatrix}
\begin{bmatrix}
\vdots \\ y_{-1} \\ y_0 \\ y_1 \\ \vdots
\end{bmatrix}
\end{gather}

For the approximate solution 
of the biinfinite Toeplitz system in~\eqref{toepsys} by the finite section 
method one considers the finite linear equations
\begin{equation}
\AAA_{N} x^{(N)} = y^{(N)}
\label{fineq}
\end{equation}
where $y^{(N)}=(y^{(N)}_{-N},\dots,y^{(N)}_{N})$ and where $\AAA_{N}$ are the
finite sections of the biinfinite Toeplitz matrix in~\eqref{toepsys}:
\begin{gather}
\AAA_{0} = (\aaa_{0})\,, \,
\AAA_{1} = 
\begin{bmatrix}
\aaa_{0}  & \aaa_{2}  & \aaa_{2} \\
\aaa_{-1} & \aaa_{1}  & \aaa_{1} \\
\aaa_{-2} & \aaa_{-1} & \aaa_{0} \\
\end{bmatrix}
\,, \, \AAA_{2} = 
\begin{bmatrix}
\aaa_{0}  & \aaa_{1}  & \aaa_{2}  & \aaa_{3}  & \aaa_{4} \\
\aaa_{-1} & \aaa_{0}  & \aaa_{1}  & \aaa_{2}  & \aaa_{3} \\
\aaa_{-2} & \aaa_{-1} & \aaa_{0}  & \aaa_{1}  & \aaa_{2} \\
\aaa_{-3} & \aaa_{-2} & \aaa_{-1} & \aaa_{0}  & \aaa_{1} \\
\aaa_{-4} & \aaa_{-3} & \aaa_{-2} & \aaa_{-1} & \aaa_{0} 
\end{bmatrix}
\, \dots \notag
\end{gather}

%Defining the orthogonal projections $\PN$ by
%\begin{equation}
%\label{defP}
%\PN x = (\dots, 0,0, x_{-N}, x_{-N+1},\dots , x_{N-1}, x_N, 0,0, \dots)
%\end{equation}
%and identifying the image of $\PN$ with the $2N+1$-dimensional space
%$\Cst^{N+1}$ we can express the equations~\eqref{toepsys} in the form
Using the the orthogonal projections defined in~\eqref{defP}
and identifying the image of $\PN$ with the $2N+1$-dimensional space
$\Cst^{N+1}$ we can express the equations~\eqref{toepsys} in the form
\begin{equation}
\PN \AAA \PN x^{(N)} = \PN y\,.
\label{}
\end{equation}
We say that the finite section method is {\em applicable} to $\AAA$, if 
beginning with some $N \in \Nst$ for each $y \in \ltZ$ the equation 
\begin{equation}
\PN \AAA \PN x^{(N)} = \PN y  
\label{applicable}
\end{equation}
has a unique solution $x^{(N)} \in \mbox{Im}\,\PN$ and as 
$N \rightarrow \infty$ the vectors $x^{(N)}$ tend to the solution of 
$\AAA x=y$. 

Now let $\aaa$ be a continuous $\DDD \times \DDD$-matrix valued function 
on $\Tst$, $\aaa \in \Csp_{\DDD \times \DDD}(\Tst)$. Let 
$\aaa_{k} \in \Cst^{\DDD \times \DDD}$ stand for the $k$-th Fourier 
coefficient of $\aaa$, then the operator $\AAA$ defined by the action
\begin{equation}
(\AAA x)_k = \sum_{j \in \Zst} \aaa_{k-j} x_j \,, \qquad x \in \ltsp_{\DDD}
\label{blocklaurent}
\end{equation}
is called a {\em block Laurent operator} with generating function $\aaa$.
A block Laurent operator can be represented by a biinfinite block Toeplitz 
matrix, where the $\aaa_{k}$ in~\eqref{toepsys} are $\DDD \times \DDD$ 
matrices. Analogous to above, for the approximate solution of biinfinite 
block Toeplitz systems we consider finite sections of the form
\begin{equation}
\PN \AAA \PN x^{(N)} = \PN y\,.
\label{}
\end{equation}

It is a well known fact and a direct consequence of the definition of 
shift-invariant systems, that the associated frame operator $S$ is a 
block Laurent operator, since $S$ commutes with the shift operator 
$U_{\DDD} (\phi_k)_{k \in  \Zst} := (\phi_{k-\DDD})_{k \in \Zst}$.
Hence the operator equation $S\gam =\gm$ can be expressed as a
biinfinite block Toeplitz system of equations,
where the blocks of the system matrix are of size $\DDD \times \DDD$.
The blocks themselves are not Toeplitz. It is exactly this 
block-Laurent structure of $S$ that is utilized in the
polyphase representation proposed in the filter bank theory~\cite{VK95}. 

Let $\SN$ be the frame operator associated with the truncated 
frame elements $\PN \gmn$ as defined in~\eqref{pnspn}. Since 
$\SN = \PN F^{\ast} F \PN = \PN S \PN$, the corresponding
system of equations $\SN \gnam= \gnm$ is just a finite section of the 
biinfinite system $S\gamma = g$
Hence it is natural to ask if we can compute approximate the duals $\gam$
by solving
\begin{equation}
\PN S \PN \gam^{(N)} = \PN g
\label{framefin}
\end{equation}
for increasing $N$, or in other words, we want to know if the finite 
section method is applicable to $S$. Moreover we are 
also interested in the rate of convergence, in case the finite section
method applies to $S$.

\section{Rate of convergence using the finite section method} 
\label{s:rate}

If the functions $\gm$ have a certain rate of decay, one would expect
that the dual functions also share this behavior.  
We consider the following three types of decay for $\gm$ and study
first how a certain decay property translates to decay properties of
the frame operator $S$.
\begin{itemize}
\item[a)] {\em Exponential decay}: $|\gm(k)| \le c \lambda^{\alpha |k|}$
      for $0 < \alpha <1$, $\lambda \in (0,1)$ and some constant 
      $c$. 
\item[b)] {\em Polynomial decay}: $|\gm(k)| \le c (1+|k|)^{-\alpha}$ for 
      $\alpha >1$ and a constant $c$. 
\item[c)] {\em Compact support}: $\gm(k) = 0$ for $|k| > s$. 
\end{itemize}

The following proposition is a simple modification of well-known
results used in connection with the construction of splines and wavelet
bases.
\begin{proposition}
\label{proppolexpdecay}
Assume that the functions $\gmn, m=0,\dots,M-1$ constitute a 
shift-invariant frame for $\ltZ$ with frame operator $S$ and duals $\gnam$.\\
(a) Exponential decay: If $|\gm(k)| \le c_1 \lambda^{\alpha_1 |k|}$
    for $0 < \alpha_1 < \alpha <1$ and $\lambda \in (0,1)$. Then there 
    exists an $\alpha_2 < \alpha$ and a constant $c_2(\alpha_2)$ such that
    \begin{equation}
    \label{exponentialmatrix}
    |S_{k,l}| \le c_2(\alpha_2) \lambda^{\alpha_2 |k-l|} \qquad
    \mbox{for} \,\,0 < \alpha_2 \le \alpha\,.
    \end{equation}
    Further it holds for the duals $\gam$ 
    \begin{equation}
    \label{dualexp}
    |\gam(k)| \le c_3(\alpha_3) \lambda^{\alpha_3 |k|}\,.
    \end{equation}
    for some $\alpha_3< \alpha$ and a constant $c_3(\alpha_3)$.\\
(b) Polynomial decay: If $|\gm(k)| \le c (1+|k|)^{-\alpha}$ for 
    $\alpha >1$, then there exists a constant $c_1$ such that
    \begin{equation}
    \label{polynomialmatrix}
    |S_{k,l}| \le c_1 (1+|k-l|)^{-\alpha}\,.
    \end{equation}
    Further there is a constant $c_2$ such that
    \begin{equation}
    \label{dualpol}
    |\gam(k)| \le c_2 (1+|k|)^{-\alpha}\,.
    \end{equation}
\end{proposition}
\begin{proof}
The exponential decay~\eqref{exponentialmatrix} and the polynomial 
decay~\eqref{polynomialmatrix} respectively, follow from Proposition~1 
in~\cite{Jaf90} by Jaffard. Proposition~2  and Proposition~3 in~\cite{Jaf90}
imply formulas~\eqref{dualexp} and~\eqref{dualpol}.
\end{proof}

We admit that~\eqref{dualexp} and~\eqref{dualpol} are not
very satisfactory from a practical point of view, since
both estimates do not provide explicit expressions for the involved constants.
We hope to address this problem in our future work.

Fortunately for compactly supported $\gm$ we can do better.
Compact support of the $\gm$ allows us to compute the transform
coefficients $\langle f, g_{m,n} \rangle$ exactly and also the matrix 
entries of $S$ can be calculated by finite sums. In filter bank design 
the analysis filters are often designed to be compactly supported 
(FIR filter banks). Note however that compact support
of the $\gm$ does not imply compact support of the duals. 
For Gabor frames B\"olcskei has shown~\cite{Boe97} that
the minimal norm duals are only compactly supported in very specific
cases. It is well known that the redundancy of frames provides some
design freedom for the dual systems, which for instance can be used
to construct dual functions with compact support~\cite{WR90,DLL95,Boe97}.
The following results allow us to obtain a good estimate for the
decay of the duals $\gamn$ for compactly supported $\gmn$.
But before we can state the results we need some preparation. 

Recall the a matrix $T$ is called $m$-banded, if $T_{k,l} = 0$ for $|k-l|>s$.
The following theorem about the decay of the inverse of a band matrix
is due to Demko, Moss, and Smith~\cite{DMS84}. 
\begin{theorem}
\label{thdemko}
Let $T$ be a positive definite, $m$-banded, bounded and boundedly
invertible matrix in $\ltsp(I)$, where ${\cal I}=\Zst, \Zst^+$ or
$\{0,1,\dots,N-1\}$. Let $B = \|T\|$ and 
$A = \|T^{-1}\|$ and set $\kap = B/A, q = \frac{\skap-1}{\skap+1}$
and $\lam = q^{\frac{1}{m}}$. Then we have
\begin{equation}
|T^{-1}_{k,l}| \le D \lamkl
\label{banddecay}
\end{equation}
where
\begin{equation}
D = \frac{1}{A} \max\{1,\frac{(1+\skap)^2}{2\kap} \}
\label{lambdadef}
\end{equation}
\end{theorem}
Further we shall need following corollary
\begin{corollary}
\label{cordecay}
Let $T$ be a matrix whose entries decay exponentially off the diagonal,
i.e.,
\begin{equation}
|T_{k,l}| \le c \lambda^{|k-l|} \qquad k,l \in \Zst
\label{matexpdecay}
\end{equation}
and let $y$ be a sequence in $\ltZ$ with $\|y\|=1$ and with compact support 
such that $y_k= 0$ for $|k| >s$. Denote $\phi = Ty$, then it holds 
\begin{equation}
|\phi_k| \le \frac{c \lam^{-s}}{1-\lam} \lam^{|k|} \,.
\label{vecxpcdecay}
\end{equation}
\end{corollary}
\begin{proof}
We have 
\begin{equation}
\phi_k = \sum_{l=-\infty}^{\infty} T_{k,l} y_{l} 
 = \sum_{l=-s}^{s} T_{k,l} y_{l}  \qquad k \in \Zst  \notag
\end{equation}
For reasons of symmetry we can restrict ourselves to the case 
$\phi_k$ with $k \ge 0$. First we consider the case $k\ge s$:
\begin{equation}
\label{case1}
|\phi_k| \le \sum_{l=-s}^{s} c \lamkl = c \lam^k \sum_{j=0}^{2s} \lam^{s-j}
 = c \lam^{k} \frac{\lam^{-s} - \lam^{s+1}}{1-\lam} \notag
\end{equation}
Now we consider the case $0 \le k < s$:
\begin{gather}
|\phi_k| \le \sum_{l=-s}^{s} c \lamkl = \sum_{l=-s}^{k-1} c \lam^{k-l}
+\sum_{l=k}^{s} c \lam^{l-k}  \notag\\
= c \lam^{k+s} \sum_{j=0}^{k+s-1} \lam^{-j} + c \sum_{j=0}^{s-k}\lam^{j} 
\notag\\
= c \lam^{k+s} \frac{1-\lam^{k+1}}{\lam^{k+s+1}-\lam^{k+s}} + 
  c\frac{1-\lam^{s-k+1}}{1-\lam}  \notag\\
= c \lam^{k} \frac{\lam^{-k+1} - \lam^{s+1}}{1-\lam} +
  c \lam^{k} \frac{\lam^{-k} - \lam^{s-2k+1}}{1-\lam}  \notag\\
= c \lam^{k} \frac{\lam^{-k+1}-\lam^{s+1}+\lam^{k}-\lam^{s-2k+1}}{1-\lam}
\notag
\end{gather}
Now we show that for $0 \le k \le s$ it holds
\begin{equation}
c \lam^{k} \frac{\lam^{-k+1}-\lam^{s+1}+\lam^{k}-\lam^{s-2k+1}}{1-\lam}
\le c \lam^{k} \frac{\lam^{-s} - \lam^{s+1}}{1-\lam} \,. \notag
\label{SNIydecay}
\end{equation}
Set $h(x) = \lam^{-x} + \lam^{x+1} -\lam-1$ for $x=0,1,\dots$ and note that
\begin{align}
                        \lam^{-(x+1)} &> \lam^{x+1}  \notag\\
           \lra \lam^{-(x+1)}(1-\lam) &> \lam^{x+1}(1-\lam) \notag\\
       \lra \lam^{-(x+1)} - \lam^{-x} &> \lam^{x+1} -\lam^{x+2}\notag \\
\lra \lam^{x+2}+\lam^{-(x+1)}-\lam -1 &> \lam^{x+1}+\lam^{-x}-\lam -1 \notag
\end{align}
whence $h(x)$ is strictly monotonically increasing. Since $h(0)=0$ and
by setting $x= s-k$ it readily follows that 
$$\lam^{s+1}-\lam^{-k+1}+\lam^{s-2k+1}-\lam^{-k} \le \lam^{s+1}- \lam^{-s}$$
which implies \eqref{SNIydecay} and hence
\begin{equation}
|\phi_k| \le c \lam^{|k|} \frac{\lam^{-s}}{1-\lam} \qquad \quad k \in \Zst\,.
\label{matvecdecay}
\end{equation}
\end{proof}

Now it is easy to estimate the decay of the dual functions $\gamn$ for 
compactly supported $\gm$.
\begin{proposition}
Let $\{\gm\}$ be a shift-invariant frame for $\ltZ$ with $\|\gm\|=1$
and assume that $g(k) = 0$ for $|k|>s$. Let $A,B$ be the frame bounds and set 
$\kappa = B/A$ and $q = \frac{\skap-1}{\skap+1}$
Then it holds
\begin{equation}
\label{dualband}
|\gam(k)| \le C \lam^{|k|} \frac{\lam^{-s}}{1-\lam}\,.
\end{equation}
where $\lam = q^{\frac{1}{2s}}$ and 
$C = \frac{1}{A} \max\{1,\frac{(1+\skap)^2}{2\kap}\}$.
\end{proposition}
\begin{proof}
Observe that $S$ is a $2s$-banded matrix. Now the result is an immediate 
consequence of Theorem~\ref{thdemko} and Corollary~\ref{cordecay}. 
\end{proof}

We are now able to answer the question about the approximation of dual 
shift-invariant frames by duals computed via the finite section method.
\begin{theorem}
\label{theorem1}
Assume that $\{\gmn\}$ is a shift-invariant frame for $\ltZ$ with 
$\|\gm\|=1$ and frame bounds $A,B$. Set $\SN = \PN S \PN$ and 
$\gnam = \SNI \PN \gm$ for $m=0,1,\dots,M-1$ and denote $\kappa=B/A$. 
Then it holds:\\
(i) The finite section method is applicable to $S$ and
  $\gnam$ converges to $\gam$ for $N \rightarrow \infty$. \\
(ii) If the $\gm$ are compactly supported with $\gm(k) = 0$ 
  for $|k|>s$ and $N>2s$, then the rate of convergence can be estimated by
  \begin{equation}
  \|\gam -\gnam \| \le  \sqrt{2} C \lam^N (\lam^{s}-\lam^{s+1})^{-3}
  \label{gogondiff}
  \end{equation}
  where $q = \frac{\skap -1}{\skap+1}$, $\lambda = q^{\frac{1}{2s}}$ and
  \begin{equation}
  \label{constantdef}
  C = \frac{1}{A} \max\{2\kappa,(1+\skap)^2 \}\,.
  \end{equation}
(iii) If $|\gm(k)| \le c_1 \lambda^{\alpha_1 |k|}$
  for $0 < \alpha_1 < \alpha <1$ and $\lambda \in (0,1)$, then there exists 
  an $\alpha_2 < \alpha$ and a constant $c_2(\alpha_2)$ such that
  \begin{equation}
  \|\gam -\gnam \| \le \kappa c_2(\alpha_2) \lambda^{\alpha_2 N}\,.
  \label{gogondiffexp}
  \end{equation}
(iv) If $|\gm(k)| \le c (1+|k|)^{-\alpha}$ for $\alpha >1$, then 
  there exists a constant $c_1$ such that
  \begin{equation}
  \|\gam -\gnam \| \le \kappa c_1 (1+|N|)^{-\alpha}\,.
  \label{gogondiffpol}
  \end{equation}
\end{theorem}

\begin{proof}
(i) Recall that $S$ is a block Laurent operator.
It follows from Theorem~9.3 in~\cite{GF74} by Gohberg and Fel'dman 
or alternatively from Theorem~4.1 in~\cite{GK94} 
by Gohberg and Kaashoek that $\SNI$ converges to $\SI$ strongly, i.e.,  
$$\|\SN x - S x\| \rightarrow 0 \quad \forall x \in \ltZ, \,\,
\text{for} \,\,N \rightarrow \infty$$ if the generating function of $S$
is continuous on the unit circle $\Tst$ and
if the operators $Q S Q$ and $Q S^* Q$ are invertible on $\ltsp(\Zst_{+})$,
where the projection $Q$ is defined by
$$Qx = (\dots,0,x_0,x_1,\dots)\,.$$
The invertibility of $QSQ$ is a direct consequence of the 
positive-definiteness of $S$, since for $x\in \mbox{Im}\,Q, \|x\|=1$, 
we have $\langle QSQ x,x \rangle = \langle S x, x \rangle > 0$,
hence $QSQ$ is invertible. In fact using the same argument we conclude 
that $\SN$ is invertible for all $N=0,1,2,\dots$.

The continuity of the generating function of $S$ follows from the fact 
that $S$ is banded, thus the generating function is a trigonometric matrix 
polynomial and therefore continuous on the unit circle $\Tst$. Hence the 
projection method is applicable to $S$ for any right hand side $\in \ltZ$
and the first part of the theorem is proved.

In order to estimate the rate of convergence for (ii), (iii) and (iv)
consider 
\begin{align}
\|\gam - \gnam\| & = \|\SI \gm - \SNI \PN \gm\| = 
\|\SI \gm - \SI S \SNI \PN \gm\| \notag \\
& \le \|\SI\| \left(\|\gm - \SN \SNI \PN \gm\| + 
   \|\SN \SNI \PN \gm - S \SNI \PN \gm\| \right) \notag \\
& \le \|\SI\| \left( \|\gm - \PN \gm\| + 
  \|(\SN - S) \SNI \PN \gm \| \right)\notag \\
& \le A^{-1} \|(\SN - S) \SNI \PN \gm \| 
\label{mainineq}
\end{align}

(ii) Recall that $S$ is a $2s$-banded matrix, i.e.,
$S_{k,l} = 0$ for $|k-l| > 2s$. Applying Theorem~\ref{thdemko} we can 
estimate the decay of the entries of $\SI$ by
\begin{equation}
\label{frameinvdecay}
|\SI_{k,l}| \le D \lamkl
\end{equation}
where $D$ and $\lam$ are as in Theorem~\ref{thdemko} with $m=2s$. 
Let $\mu^{(N)}_0,\dots,\mu^{(N)}_{N-1}$ be the eigenvalues of $\SN$.
Since $\SN$ is Hermitian positive definite, it follows from Cauchy's 
Interlace Theorem~\cite{HJ94} that the eigenvalues 
$\mu^{(N-1)}_0,\dots,\mu^{(N-1)}_{N-1}$ of $S_{N-1}$ satisfy
$$\mu^{(N)}_0 \le \mu^{(N-1)}_0 \le \mu^{(N)}_1 \le \dots \le \mu^{(N-1)}_{N-2}
\le \mu^{(N)}_{N-1}\,,$$
which implies 
$$ \cond(\SN) \le \cond(S)\,.$$
Hence we can use the same $\lambda$ and the same constant $D$ in order to 
bound the decay of the entries of $\SNI$ independently of $N$ and obtain
\begin{equation}
\label{frameinvdecay1}
|(\SNI)_{k,l}| \le D \lamkl \,.
\end{equation}

Set $\phinm = \SNI \PN \gm$ then it follows from Lemma~\ref{cordecay} that
\begin{equation}
|(\phinm)_{k}| \le D \lam^{|k|} \frac{\lam^{-s}}{1-\lam} \,\,,\qquad \quad 
k \in \Zst\,.
\label{matvecdecay1}
\end{equation}
In the next step we estimate $\|(\SN-S)\SNI \PN \gm\|$. Write $\en = (\SN - S) 
\phin$ and observe that $\|\en\|^2 =\sum_{k=-\infty}^{\infty}|(\en)_{k}|^2
=2\sum_{k=0}^{\infty}|(\en)_{k}|^2$ (note that $(\en)_{0} = 0$). For
$k=0,1.,\dots$ consider
\begin{equation}
(\en)_{k} = \sum_{l=-\infty}^{\infty} (\SN - S)_{k,l} (\phinm)_{l} =
\sum_{l=k-2s}^{k+2s} (\SN - S)_{k,l} (\phinm)_{l}\notag
%\label{}
\end{equation}
Clearly $(\en)_{k} = 0$ for $k \le N$, hence we only have to consider the
case $k > N$:
\begin{align}
|(\en)_{k}| & = |\sum_{l=k-2s}^{k+2s} (\SN - S)_{k,l} (\phinm)_{l} | \le
 B \frac{D \lam^{-s}}{1-\lam} \sum_{l=k-2s}^{k+2s} \lam^{|l|} \notag \\
& \le \frac{B D \lam^{-s}}{1-\lam} \frac{\lam^{k-2s}-\lam^{k+2s}}{1-\lam}
\le \frac{B D}{(1-\lam)^2} \lam^{k-3s}\,, \notag 
\end{align}
whence
$$|(\en)_{k}|^2 \le \left(\frac{BD}{(1-\lam)^2}\right)^2 \lam^{2k-6s}\,.$$
Due to $(\SN-S)_{k,l} = 0$ for $|k-l| > 2s$ and since
$(\phinm)_{k} = 0$ for $|k|>N$ we obtain
$\sum_{k=0}^{\infty}|(\en)_{k}|^2 = \sum_{k=N+1}^{N+2s}|(\en)_{k}|^2$
and therefore
\begin{gather}
\|\en\|^2 = 2\sum_{k=N+1}^{N+2s} |(\en)_{k}|^2 \le 
2 \frac{B^2 D^2}{(1-\lam)^4} \lam^{-6s}\sum_{k=N+1}^{N+2s} \lam^{2k} \notag\\
\le \frac{2 B^2 D^2}{(1-\lam)^4} \lam^{2N} 
    \frac{\lam^{2-6s}-\lam^{2-2s}}{1-\lam^2} 
\le \left(\sqrt{2} BD \lam^{N} \frac{\lam^{-3s}}{(1-\lam)^3}\right)^2 \notag
\end{gather}
which implies
\begin{equation}
\|\en\| \le \sqrt{2} BD \lam^N (\lam^{s}-\lam^{s+1})^{-3}
\label{errorek}
\end{equation}
Combining~\eqref{errorek} with~\eqref{mainineq} 
yields the bound~\eqref{gogondiff}.

The statements in (iii) and (iv) follow by applying 
Proposition~\ref{proppolexpdecay} to~\eqref{mainineq}.
\end{proof}

\noindent
{\bf Remarks:}
\vspace*{-3mm}
\begin{itemize}
\setlength{\itemsep}{-0.5ex}
\setlength{\parsep}{-0.5ex}
\item[(i)] In spite of these results, good a priori estimates of the frame 
bounds as given e.g.~in~\cite{Dau92,Jan97a} play an important role.
\item[(ii)] A different approach to approximate the inverse of the frame 
operator by finite dimensional methods is due to Casazza and Christensen.
In~\cite{CC98} they consider a more general setting, however at the 
cost that their estimates are less explicit and therefore less useful
in applications.
\end{itemize}

\section{Approximation of tight frames} \label{s:tight}

An important special case of frames are so-called {\em tight frames}, which 
are defined in terms of equality of the upper and lower frame bound, 
i.e., $A=B$. Besides the nice stability property $cond(S)=1$, tight frames are
distinguished by a simple expression for the duals, which are just scaled
versions of the original frame elements, more precisely one has
$\gamn = \frac{1}{A} \gmn.$
Due to this close relationship to orthogonal bases the construction
of tight frames is of great interest in theory and in applications.
Daubechies~\cite{Dau90} has shown that for a given frame, say $\{\phi_k\}$
the set of functions $\{S^{-\frac{1}{2}}\psi_k\}$ constitutes a
tight frame. Recall that $S$ is positive definite, hence $\SR$ and $\SRI$ 
are uniquely determined. 

Recall that $S^{\frac{1}{2}}$ commutes with the same operators
that commute with $S$. For shift-invariant systems this implies that
the tight frame generated by $\phimn = S^{-\frac{1}{2}} \gmn$ 
also constitutes a shift-invariant system. We want to approximate the 
elements of the tight frame $\{\phimn\} $ by the solutions of the finite 
sections
$$(\PN S \PN)^{\frac{1}{2}} \phinm = \PN \gm\,.$$
One can easily check that the entries of $\SR$ and $\SRI$ decay 
exponentially off the diagonal, see e.g.~\cite{JM89,Lai94}.
The applicability of the finite section method follows now
from Theorem~5.1 in~\cite{GK94}. Along the same lines of the proof of 
Theorem~\ref{theorem1} one can show that 
$\phinm := (\PN S \PN)^{-\frac{1}{2}} \PN \gm$ converges exponentially
fast to $\phim$. Since we do not have explicit expressions for the
decay of the entries of $\SR$ and $\SRI$ we only can give following 
qualitative result for the approximation error.
\begin{theorem} 
\label{theorem3}
Let $\{\gmn\}$ be a shift-invariant frame for $\ltZ$ with 
frame operator $S$ and frame bounds $A,B$. Assume that $\|g\|=1$ and 
set $\phim = \SRI \gm$ and $\phim^{(N)} = (\PN S \PN)^{-\frac{1}{2}} \PN g$.
Then the finite section method is applicable to $\SR$. If
moreover $g$ is compactly supported with $g(k) = 0$ for $|k|>s$ then
$$\|\phim -\phim^{(N)} \| \le K \beta^{N} \qquad N=0,1,\dots$$
for $\beta \subset (0,1)$ and some constant $K>0$, where $\beta$ and $K$
depend on the frame bounds $A,B$ and the support length $s$, but
are independent of $N$.
\end{theorem}

\section{Rate of convergence by ``canonical'' approximation using a 
periodic model} 
\label{s:periodic}

We have seen in Section~\ref{s:rate} that the functions $\gnamn$ which are
the duals to the truncated functions $\PN \gmn$ converge exponentially fast 
to the duals $\gamn$ of the infinite-dimensional problem.
So from a theoretical point of view we know that we can compute a good
approximation to $\gam$  by solving the finite system $\SN \gnam = \gnm$.
From a numerical point of view, however, we are interested in solving
$\SN \gnam = \gnm$ as efficient as possible. 

Unfortunately due to the
truncation of $S$ we loose some structural properties of the problem. 
For instance the inverse of $S$ or the product $S^2$ still have
block Toeplitz structure, whereas $\SNI$ and $\SN^2$ are no longer
of Toeplitz type. Certainly one could apply one of the
fast solvers for Toeplitz systems, but we can also try to find a
more canonical way to design a finite dimensional model, which allows
to preserve the perfect symbol calculus of Laurent operators. 
This is one motivation to use a periodic model for
the design of numerical algorithms. 

Despite the fact that the periodic model is widely used in connection with
shift-invariant systems, almost no attempt has been made to  
investigate the relation between the duals of the periodic finite model
and the duals of the infinite dimensional problem.

One exception is~\cite{Jan97a} where Janssen shows 
that if the $\gmn$ constitute a Gabor frame for $\ltZ$ with frame
bounds $A,B$ and if $\gmn \in \liZ$, then the periodized functions
$\gmn^{per}(k)= \sum_{j=-\infty}^{\infty}\gmn(k-jL)$ constitute 
a frame for $\Cst^L$ with frame bounds $A,B$.
Furthermore the duals for the periodized frame can be obtained by 
periodizing the duals of the original frame elements. The drawback of
this otherwise very appealing result is that the length of the period 
$L$ has to be the least common multiple of the time shift parameter $a$ and 
the number of frequency channels $M$. This condition can be
a serious restriction in practice.  For instance take $a=2, M=3$, then 
the periodized functions are of length 6, which may be much too short for 
many applications.

In this section we will clarify the question under which conditions the 
duals computed within the periodic model converge to the duals
$\gamn$ when the length of the period tends to infinity.

Let the functions $\gm$ be compactly supported, such that $\gm(k) = 0$
for $|k|>s$. We construct our periodic model, as usual, by
extending the truncated  functions $\PN \gm$ periodically
beyond the interval $[-N,N]$ for $N >s$ by setting
\begin{equation}
\label{cgm}
\cgm(k+lN) := \gm(k) \qquad \text{for}\,\, k=-N,\dots,N,\,\,, l \in \Zst\,. 
\end{equation}
The other elements of our periodic finite shift-invariant system are now 
given by
\begin{equation}
\label{cgmn}
\cgmn(k) = \cgm(k-na)
\end{equation}
where the shift, according to our model, is understood as circulant shift.
The periodicity of the $\cgmn$ implies that the corresponding frame operator 
is a block circulant matrix, rather than merely a block Toeplitz matrix.
Observe that the entries of $\CN$ deviate from the entries of $\SN$ only 
in the lower right and upper left corner of the matrix.

In other words we have replaced the block Toeplitz matrix $\SN$ in
the system $\SN \gnam = \gnm$ by a block circulant matrix.
However replacing a block Toeplitz system by a block circulant matrix
in a linear system of equations will in general significantly affect 
the solution. 

Moreover it is a priori not clear under which conditions invertibility 
of $\SN$ implies invertibility of $\CN$. For instance take following 
Hermitian positive definite Toeplitz matrix $T$ and its 
``circulant completion'' ${^{\footnotesize P}T}$:
\begin{gather}
T = 
\begin{bmatrix}
2  & -1 & 0  & 0  \\
-1 &  2 & -1 & 0  \\
0  & -1 &  2 & -1 \\
0  &  0 & -1 & 2  \\
\end{bmatrix}
\,, \qquad  {^{\footnotesize P}T} =
\begin{bmatrix}
2  & -1 & 0  & -1 \\
-1 &  2 & -1 & 0  \\
0  & -1 &  2 & -1 \\
-1 &  0 & -1 & 2  \\
\end{bmatrix}\notag
\end{gather}
then $T$ is invertible, whereas ${^{\footnotesize P}T}$ is obviously
not invertible.

The following theorem justifies the usage of the circulant model for 
the approximation of dual shift-invariant systems. Implicitly the theorem 
also provides an estimate of the error we make when we 
approximate the solution of a block Toeplitz system by solving the
related block circulant system.
\begin{theorem}
\label{theorem2}
Let $\{\gmn\}$ be a shift-invariant frame for $\ltZ$ with 
frame bounds $A,B$. Assume that $\|\gm\|=1$ and that the $\gm$ are compactly 
supported with $\gm(k) = 0$ for $|k|>s$. Let $N>2s$ and let $\cgmn$
be the periodized functions as defined in~\eqref{cgm} and~\eqref{cgmn} with 
frame operator $\CN$. Then $\CN$ is invertible and $\cnam = \CNI \cgm$ 
converges to $\gam$ for $N \rightarrow \infty$. If $N>3s$ then the rate 
of convergence can be estimated by
\begin{equation}
\|\PN \gam -\cnam \| \le  3 \sqrt{2} C \lam^N (\lam^{s}-\lam^{s+1})^{-3}
\label{gogondiff2}
\end{equation}
where $\kappa = \frac{B}{A}, q = \frac{\skap -1}{\skap+1}$,
$\lambda = q^{\frac{1}{2s}}$ and
\begin{equation}
\label{constantdef1}
C = \frac{1}{A} \max\{2\kappa,(1+\skap)^2 \}
\end{equation}
\end{theorem}

\begin{proof}
First we show that invertibility of $S$ implies invertibility of $\CN$.
The results about the spectrum of block circulant matrices in the proof 
of Theorem~3 in~\cite{ZTA97} in combination with Corollary~1 
in~\cite{ZTA97} imply that $$\sigma(\CN) \subseteq \sigma(S)\,.$$
$\sigma(S)$ is bounded away from zero since $A>0$, hence
it follows that $\CN$ is invertible.

%From Theorem~\ref{theorem1} we know that $\gn$ approximates $\go$ 
%for $N \rightarrow \infty$, where the rate of convergence is given by 
%\begin{equation}
%\|\go -\gon \| \le  \sqrt{2} C \lam^N (\lam^{s}-\lam^{s+1})^{-3}\,.
%\end{equation}
%We investigate how well $\gon$ is approximated by $\con$
%using the circulant model.
To prove~\eqref{gogondiff2} we write
\begin{equation}
\|\PN \gam - \cnam \| \le 
%\|\PN \gam - \PN \gnam \|+\|\PN \gnam - \cnam\| \le
 \|\gam - \gnam \|  + \|\PN \gnam -\cnam\|\,.
\label{circdiff}
\end{equation}
We already have estimated $\|\gam - \gnam \|$ in~\eqref{gogondiff} of
Theorem~\ref{theorem1}, so it remains to estimate $\|\PN \gnam -\cnam\|$. 
Proceeding analogously to 
the calculation preceding equation~\eqref{mainineq} we can write
\begin{equation}
\|\gnam - \cnam\| \le A^{-1} \|(\CN - \SN) \CNI \PN \gm\|\,.
\label{mainineq2}
\end{equation}
Note that $\CN$ is three-band
matrix, with one band centered at the main diagonal, and two other bands of 
width $2s$ located at the lower left and upper right corner of the matrix.
It follows from Proposition~5.1 in~\cite{DMS84} that the entries of 
$\CNI$ decay exponentially off the diagonal and off the lower right and
upper left corner. More precisely 
\begin{equation}
|(\CNI)_{k,l}| \le 
\begin{cases}
D \lam^{|k-l|} & \text{if $0 \le |k-l| \le N$} \\
D \lam^{2N+1-|k-l|} & \text{if $N+1 \le |k-l| \le 2N$} 
\end{cases}
\label{Scond}
\end{equation}
where $D$ and $\lam$ are as in Theorem~\ref{thdemko} with $m=2s$.

Before we show that $\phinm: = \CNI \cgm$ also decays exponentially,
observe that due to the construction of $\CN$ the matrix
$\EN: = \CN - \SN$ is a sparse matrix having non-zero entries only
in the upper right and the lower left corner, as illustrated in~\eqref{CNSNmat}
\begin{gather}
\qquad \qquad \begin{matrix}
& & & & & \text{\tiny{$N\!-\!2s+1$}} & & & \text{\tiny{$N$}} \notag 
\vspace*{-6mm}
\end{matrix} \\
\EN = 
\begin{bmatrix}
   & & &           & \xx       & \dots  & \xx       \\
   & & &           &         & \ddots & \vdots  \\
   & & &           &         &        & \xx       \\
   & & & \text{\Large{$0$}}  &         &        &         \\
 \xx & & &           &         &        &         \\
\vdots & \ddots & & &        &        &         \\
\xx & \dots & \xx &    &         &        &         \,\,\,.
\label{CNSNmat}
\end{bmatrix} 
%\vskip-20mm
%\begin{matrix}
%\text{\tiny{$-N$}} & & \text{\tiny{$-N-1+2s$}} & & & & & &  \notag 
%\end{matrix} 
\end{gather}
Denoting $\en = \EN \CNI \cgm$, it readily follows from the sparsity
of $\EN$ that only the first $2s$ and the last $2s$ components of
$e_N$ are non-zero, hence only the decay of the first $2s$ and last $2s$ 
components of $\phinm = \CNI \cgm$ is of interest.

Since $\CN$ is a three-band matrix, the decay behavior of the entries
of $\CN^{-1}$ is a little bit more complicated than for ordinary band 
matrices. We split the analysis of the decay of the entries of $\phinm$ into 
two steps by considering first the entries $(\phinm)_{k}$ for 
$k = N-2s,\dots,N-m-1$ and then those for $k=N-m,\dots,N-1$.

\noindent
Case $k=N-2s+1,\dots,N-s$: \\
$$(\phinm)_{k} = \sum_{l=-N}^{N} (\CNI)_{k,l} (\cgm)_l = 
\sum_{l=-s}^{s} (\CNI)_{k,l}(\cgm)_l$$ due to the
compact support of $\gm$. Hence 
\begin{equation}
|(\phinm)_{k}| \le |\sum_{l=-s}^{s}(\CNI)_{k,l}|
\le D \sum_{l=-s}^{s} \lam^{|k-l|}  \,.
\label{}
\end{equation}
The assumption $N > 3s$ implies that $k \ge l$ hence
\begin{align}
D \sum_{l=-s}^{s} \lam^{|k-l|} = D \sum_{l=-s}^{s} \lam^{k-l} 
= D \lam^{k} \frac{\lam^{-s}-\lam^{s+1}}{1-\lam}
\le \frac{D}{1-\lam}\lam^{k-s}
\end{align}

\noindent
Case $k=N-m,\dots,N-1$: \\
\begin{align}
|(\phinm)_{k}| & \le |\sum_{l=-s}^{s}(\CNI)_{k,l}| 
\le \sum_{l=-s}^{k-N-1}|(\CNI)_{k,l}| +
\sum_{l=k-N}^{s}|(\CNI)_{k,l}| \notag \\
& \le D \sum_{l=-s}^{k-N-1} \lam^{2N+1-k+l} 
+ D \sum_{l=k-N}^{s} \lam^{k-l} \notag \\
& = D \lam^{2N+1-k-s} \sum_{l=0}^{k-N+s-1} \lam^{l}
+ \frac{D}{1-\lam} \lam^{k} \frac{\lam^{-s}-\lam^{-(k-N-1)}}{1-\lam}\notag \\
& =\frac{D}{1-\lam}
\left(\lam^{2N+1-k-s}-\lam^{-N+1}+\lam^{k-s}-\lam^{N+1}\right) \notag \\
& \le \frac{D}{1-\lam}(\lam^{k-s}+\lam^{2N+1-k-s}) \le 
\frac{2D}{1-\lam}\lam^{k-s} 
\end{align}
Note that 
$\|\en\|^2 = \sum_{k=-N}^{N}|(\en)_{k}|^2 = 2\sum_{k=N-2s+1}^{N}|(\en)_{k}|^2$.
Due to the special sparsity structure of $\EN$ as illustrated 
in~\eqref{CNSNmat} we can write
\begin{align}
|(\en)_{k}| & \le \sum_{l=2N-2s+1+k}^{N} |(\EN)_{k,l}| |\phi_l|
\le \frac{2BD}{1-\lam} \lam^{-s} \sum_{l=2N-2s+1+k}^{N} \lam^l \notag \\
& \le \frac{2BD}{(1-\lam)^2} \lam^{-s} (\lam^{2N-2s+1+k} - \lam{N+1} 
\le \frac{2BD}{(1-\lam)^2} \lam^{k-3s+2N} \notag 
\end{align}
It follows
$$|(\en)_{k}|^2 \le \left(\frac{2BD}{(1-\lam)^2}\right)^2 \lam^{2k+4N-6s}$$
and  
\begin{gather}
\|\en\|^2  =  2 \sum_{k=-N}^{-N+2s-1} |(\en)_{k}|^2 
\le 2 \frac{(2BD)^2}{(1-\lam)^4} \sum_{k=-N}^{-N+2s-1} 
\lam^{2k-6s+4N} \notag \\
\le 2 \frac{(2BD)^2}{(1-\lam)^4} 
\frac{\lam^{2N-6s} - \lam^{2N-2s}}{1-\lam} 
\le \left(2 \sqrt{2} \frac{BD}{(1-\lam)^3} \lam^{N-3s}\right)^2 \notag
\end{gather}
and therefore
$$
\|\en\| \le 2 \sqrt{2} \frac{BD \lam^{-3s}}{(1-\lam)^3} \lam^{N}
= 2 \sqrt{2} BD (\lam^{s}-\lam^{s+1})^{-3} \lam^{N}
$$
which together with~\eqref{mainineq2} yields the desired 
estimate~\eqref{gogondiff2}.
\end{proof}

The corresponding variant of Theorem~\ref{theorem3} for the approximation
of a tight frame via the periodic model is left to the reader.

In fact, the periodic extension of finite signals is more than just a 
simple and convenient way to hanlde boundary problems. It is in some sense 
the ``canonical'' way to set up a discrete model, when translations
come into play, since it allows to preserve important mathematical 
properties of the infinite dimensional problem.

More generally speaking one preserves the underlying group structure 
of the original problem. By doing so, we can immediately apply the
abstract results derived for shift-invariant systems~\cite{FK98,Gro98,TA98}. 
In fact, sometimes, (e.g.~in the case of
Gabor analysis~\cite{Str98}) we gain even more structural properties compared
to the infinite dimensional problem, due to the finiteness of the 
underlying abelian group.

\section*{Acknowledgement}
I want to thank B. Silbermann and G. Zimmermann for their valuable
suggestions.

The author has been partially supported by Schr\"odinger scholarship
J01388-MAT of the Austrian Science foundation FWF.

\bibliographystyle{plain}

\end{document}